
\input amstex
\batchmode
\input amssym.def
\documentstyle{amsppt}
\topmatter
\title
cuspidal hypergeometric functions
\endtitle
\author
Eric M. Opdam
\endauthor
\dedicatory
Dedicated to Richard Askey, on the occasion of his 65th birthday
\enddedicatory
\affil
Department of Mathematics, University of Leiden
\endaffil
\address
P.O. Box 9512, 2300 RA Leiden, The Netherlands
\endaddress
\date
October 1996
\enddate
\endsubjclass
\thanks I thank Patrick Delorme and Gerrit Heckman and for many helpful conversations and comments
\endthanks
\abstract 
We prove the Plancherel formula for hypergeometric functions 
associated to a root system in the situation when the root multiplicities 
are negative (but close to 0). As a result we obtain a 
classification of the hypergeometric functions that are square 
integrable, and we find a closed formula for their square norm 
as a function of the root multiplicities.
\endabstract
\endtopmatter
\document
\baselineskip=15pt

\font\title=cmbx10 scaled\magstep1

\def \a {\alpha}
\def \b {\beta}
\def \g {\gamma}
\def \G {\Gamma}

\def \e {\epsilon}

\def \k {\kappa}
\def \l {\lambda}

\def \cci {C_c^\infty}

\def \fa {\frak a}
\def \fh {\frak h}

\def \Z {\bold Z}

\def \R {\bold R}

\def \F {\Cal F}

\def \K {\Cal K}

\beginsection{1. Introduction}

This section serves to explain the main results of this paper 
(see the last two paragraphs of this section), and to 
review some basic facts and notations of the theory of 
hypergeometric functions associated with root systems.
 
Let $\frak a$ be a Euclidean space of dimension $n$ and
$R\subset\frak a^*$ (the dual of $\frak a$) be an
integral, irreducible, reduced root system which 
spans $\fa^*$. 
Its Weyl group is denoted by $W$. Let $\K$ denote the 
linear space of multiplicity functions, i.e. the
space of $W$-invariant real functions on $R$. We 
fix a set of positive roots $R_+$. 
If $ k\in\K$ we consider the following differential 
operator on $\frak a$:
$$
L(k)=\sum_{i=1}^n\partial(X_i)^2-
\sum_{\alpha\in R_+}k_\alpha(1+e^\alpha)(1-e^\alpha)^{-1}
\partial(X_\alpha).
\tag1.1
$$
Here we use the convention to write $\partial(p)$ (with 
$p$ an element of the symmetric algebra $Sym(\frak a)$ 
of $\frak a$) 
for the constant coefficient differential operator 
on $\frak a$ that corresponds with $p$. Moreover, 
we have chosen an orthonormal basis $\{X_i\}$ of $\frak a$, 
and $X_\a\in\frak a$ is the element such that $\l(X_\a)=(\l,\a)
\,\forall\l\in{\frak a}^*$.
This remarkable operator has been studied intensely over the  
years ([6], [11], [12], [8], [16], [7], [17], [4], [10]), for 
some very good reasons. It arises "in nature" as the 
(radial part of the) Laplace-Beltrami operator of noncompact 
Riemannian symmetric spaces for special choices of the 
parameter $k$. It was realized that this operator plays 
a key role in the understanding of the Macdonald constant 
term conjectures. Finally this operator attracted some 
attention of mathematical physicists ([14], [15]) because  
of the relation of $L(k)$ with the quantum Calogero-Moser 
system. More precisely, define 
$$
\rho(k)=\rho(R_+,k)={1\over 2}\sum_{\a\in R_+} k_\a\a\
\tag1.2
$$
and 
$$
\delta(k;x)= \prod_{\alpha\in R_+}|2\sinh(\frac{1}{2}\a(x))|^{2k_\a}, 
\tag1.3
$$
then the operator $S(k)=\sqrt{\delta(k;x)}(L(k)+(\rho(k),\rho(k)))
\sqrt{\delta(k;x)}^{-1}$
takes the form 
$$
S(k)=\sum_{i=1}^n\partial(X_i)^2-
\frac{1}{4}\sum_{\alpha\in R_+}\frac{(\a,\a)k_\a(k_\a-1)}
{\sinh^2(\frac{1}{2}\a(x))}
\tag1.4
$$
This is the Schr\" odinger operator of the Calogero-Moser system.
From this last formula it is clear that $L(k)$ is a formally 
symmetric operator on the space $\cci(\fa)^W$, if we equip this 
space with the Hermitian inner product 
$$
(\phi,\psi)=\int_{\fa}\phi(x)\overline{\psi(x)}\delta(k;x)dx
\tag1.5
$$
where $dx$ denotes the Lebesgue measure on $\fa$. 
Morever, we fix the length of the roots  
in such a way that a fundamental domain of the lattice 
$Q^\vee$ (here $Q^\vee$ is the coroot lattice in $\fa$, 
spanned over $\Z$  
by the coroots $\a^\vee=\frac{2X_\a}{(\a,\a)}$ ($\a\in R$))
has volume equal to 1 (this turns out to be the most natural 
normalization as we shall see later).
In this paper we shall discuss the spectral problem for 
$L(k)$, in other words we try to decompose the $L_2$ space 
of $W$-invariant functions associated with the measure in (1.5) 
as a direct integral of eigenspaces for $L(k)$. 
Formula (1.5) makes sense only as long as $\delta(k;x)$ is a 
locally integrable function on $\fa$, and this imposes a 
condition on $k\in \K$. We shall restrict ourselves exclusively 
to this situation where $\delta(k;x)$ is indeed locally 
integrable. Obviously this condition is 
satisfied if $k_\a>0\,\forall\a\in R$, and in this situation 
the spectral problem for $L(k)$ was solved completely in [17].
In this situation the spectrum turns out to be completely 
continuous, and the spectral decomposition is very similar 
to the decomposition of the space $L^2(X)$ in irreducible 
(spherical) representations for $G$ if $X=G/K$ is a 
noncompact Riemannian symmetric space (cf. [6], [11]). The 
case that we are about to study here in this paper is the one 
in which $k_\a<0\,\forall\a\in R$, but satisfying condition 
(1.6) of Proposition 1.1 below so as to 
assure that $\delta(k;x)$ is locally 
integrable. It is elementary that if $k_\a<0\,\forall\a\in R$ and 
$\delta(k;x)$ is locally integrable then $\delta(k;x)$ is even integrable.
\proclaim{Proposition 1.1} Let  
$\b\in R$ be the highest short root of $R$. If 
$$
k_\a<0\,\forall\a\in R,\text{ and }\rho(k)(\b^\vee)+k_\b+1>0
\tag1.6
$$
then $\delta(k;x)$ is integrable.
\endproclaim
\demo{Proof} In ([2], Corollary 2.2) the following sufficient 
condition was given: $\delta(k;x)$ is integrable if 
$$
k_\a<0\,\forall\a\in R,\text{ and }\sum_{i=1}^m k_ih_i+1>0
\tag1.7
$$
Here $h_i=\#(R_i)/n$ if the $R_i$ ($i=1,\dots,m$) denote the $W$-orbits 
in $R$ (so $m=1$ or $2$). If $m=1$ then this is equivalent to (1.6) 
because the height of the highest coroot equals $h_1-1$, as is well 
known. By the same formula it is enough to check in the case $m=2$ 
that if 
$\sigma$ denotes half of the sum of the short positive roots, one has 
$\sigma(\b^\vee)=h_1-1$ if $R_1$ is the set of short roots. If $R_1$ is 
irreducible, this follows again from the above formula applied to $R_1$. 
If $R_1$ is reducible, then $R=B_n$ and here one verifies 
directly that $\sigma(\b^\vee)=h_1-1=1$.
\enddemo

In the situation of Proposition 1.1 we shall give the spectral 
decomposition of the operator (1.1) almost completely explicitly, 
and as we shall see 
this involves lower dimensional spectral families. We can handle 
the situation by a combination of the techniques of [9], 
where we dealt with the so called "Yang system" in the attractive case, 
and of [17]. In fact it turns out that one has to slightly adapt 
the arguments of these papers at two critical points, and after that 
most of the arguments of these papers will go through. It is not our 
intention to repeat all the arguments of these papers here so as to 
make this paper self-contained, because that would certainly result in 
a very lengthy story. Instead, we shall just indicate the two critical 
steps mentioned above, and after that freely use the statements of [9] and 
[17]. In any case, the results are quite nice. We obtain a complete 
classification of the square integrable eigenfunctions of $L(k)$ and we 
can, up to some rational constant, evaluate their square norms 
explicitly as a function of the multiplicities $k_\a$. Moreover, 
it seems likely that similar methods can be applied in even greater 
generality. For example, we have excluded the case $R=BC_n$ here 
because that really seems to call for more delicate arguments of 
a combinatorial nature. But that case is very interesting in its own right
and it should be investigated seperately. Also it seems likely that some 
form of the methods presented here should be applicable to the 
harmonic analysis of the so called double affine Hecke algebra introduced 
by Cherednik. In this setting one should obtain the Aomoto conjecture, 
recently proved by Macdonald [13], as the square norm computation 
of the simplest discrete eigenfunction.
 
In section 2 we shall gather together some necessary notation and 
"basic" results about $L(k)$ and its eigenfunctions.   
In section 3 we shall discuss the process of taking residues as 
it was formulated in [9] and explain how this applies to the 
situation in this paper. In section 4 we prove the Paley-Wiener 
theorem and finally in section 5 we shall put things together and 
prove the main formula, the Plancherel formula.

\beginsection{2. On the eigenfunctions of $L(k)$}

To describe the eigenfunctions of $L(k)$ we start by  
recalling one of the most important features of $L(k)$, 
namely that $L(k)$ is an element of a commutative algebra 
of $W$-invariant differential operators which is naturally 
isomorphic to $Sym(\fh)^W$ (here $\fh$ denotes the 
complexification of $\fa$). This fact can be understood in 
an elementary way using following beautiful observation 
about $L(k)$ due 
to Cherednik ([4], also see [5], [7]). Define the so-called 
Cherednik operator 
$$
D_\xi(k)=\partial_\xi+\sum_{\alpha\in R_+}k_\alpha\alpha(\xi)
       {1\over 1-e^{-\alpha}}(1-r_\alpha)-\rho(k)(\xi)
\tag2.1      
$$
(Here $r_\a$ denotes the reflection in the root $\a$, 
and $\xi$ denotes an arbitrary element of $\fh$.)
The following facts about the operators $D_\xi$ were discussed 
systematically in [17]. First of all, 
these operators map $\cci(\fa)$ to itself (despite the 
apparent singularity in (2.1)). Next, the operators mutually 
commute, and generate an operator algebra isomorphic to $Sym(\fh)$
where the isomorphism $Sym(\fh)\ni p\to p(D)(k)$ is determined by 
the assignment $\fh\ni \xi\to D_\xi(k)$. Moreover, the 
$W$-invariant operators of the form $p(D)(k)$ are exactly those 
with $p\in Sym(\fh)^W$. Now we can restrict such $W$-invariant 
operators $p(D)(k)$ to $\cci(\fa)^W$ and the resulting operator 
will be a {\sl {differential}} operator mapping  $\cci(\fa)^W$
to itself. This differential operator will be denoted by $D_p(k)$.
Finally $L(k)$ turns out to be an operator of this form:
$$
D_{{\sum X_i}^2}(k)=L(k)+(\rho(k),\rho(k)).
\tag2.2
$$

Let us now turn to the eigenfunction problem for the operators 
$D_p(k)$. For these facts we refer to [17] and the 
references therein. The simplest eigenfunctions are the asymptotically 
free eigenfunctions $\Phi(\l,k;x)$ on $\fa_-$. Indeed, if 
we assume that $\Phi(\l,k;x)$ has an asymptotic expansion
on $\fa_-$ of the form ($Q$ denotes the root lattice, 
i.e. the lattice spanned by the 
roots (over $\Z$)).
$$
\Phi(\l, k;x) =
e^{(\l+ \rho (k))(x)} \sum_{\k\in Q_+} \Delta_\k
(\l,k) e^{\k(x)}
\tag2.3
$$ 
with 
$$
\Delta_0 (\l,k) = 1
\tag2.4
$$
then all the coefficients are determined by recurrence 
by the single equation 
$$
L(k)\Phi(\l, k;x)=(\l+\rho(k),\l-\rho(k))\Phi(\l, k;x).
\tag2.5
$$
It is easy to see that $\Phi(\l, k;x)$ automatically is a 
simultaneous eigenfunction for all the $D_p(k)$:
$$
D_p(k)\Phi(\l, k;x)=p(\l)\Phi(\l, k;x).
\tag2.6
$$
\proclaim{Lemma 2.1} Let $\k\in Q$ but $\k\not=0$. 
Then $\k^\vee\in CH(R^\vee)$, the convex hull of the coroots.
\endproclaim 
\demo{Proof} We may assume that $\k\in Q^+$ and argue by 
induction on the height of $\k$. If the height of $\k$ equals 1 
the statement is obviously true. Next note that if 
the statement holds for two nonzero elements $\k_1$ and $\k_2$ in $Q^+$ and 
$(\k_1,\k_2)\geq 0$ then the statement also holds for $\k_1+\k_2$. If $\k
\in Q^+$ and $\k\not= 0$ then there exists a simple root $\a$ such that 
$(\k,\a)>0$. If $(\k-\a,\a)\geq 0$ then we are done by the previous 
observation 
applied to the decomposition $\k=\a+(\k-\a)$. However, if $(\k-\a,\a)<0$
then $\k=r_\a(\k-\a)$ and we are done by the $W$-invariance of the problem.
\enddemo

We shall need the following facts about the asymptotically free
eigenfunctions:
\proclaim{Proposition 2.2} 
\roster
\item
$\Delta_\k (\l,k)$ is a rational function, with poles at hyperplanes
of the form $\l(\k^\vee ) + 1 = 0$ $(\k \in Q_+ \backslash \{ 0 \})$ only.
\item
Let $x_0 \in \fa_+$, and $\e>0$. Then $\exists K_{x_0,\e} \in \R_+$ such that
$\forall \l$ with $\text{Re}(\l(\a^\vee))>\e-1
\ \forall \a\in R_+$ 
the following holds: 
$\forall \k\in Q_+$:
$$
\vert \Delta_\k (\l,k) \vert \leq K_{x_0,\e}\text{e}^{\k(x_0)}
$$
In particular, $\Phi(\l, k;x)$ is analytic in $\l$ if 
$\text{Re}(\l(\a^\vee))>-1$ and $x\in \fa_-$.
\endroster
\endproclaim
\demo{Proof} See [17], Theorem 6.3. We have slightly relaxed  
the usual condition on $\l$ in \therosteritem2, and if one checks the 
proof of Gangolli in the classical context (cf. [11], Ch. IV, 
Lemma 5.6) it is clear that this is justified by Lemma 2.1.
\enddemo
The next step is the introduction of the 
hypergeometric function for $R$, which is 
the unique (up to normalization) simultaneous eigenfunction of the $D_p(k)$ 
that extends holomorphically to a tubular neighbourhood of $\fa$ in $\fh$. We need to define 
Harish-Chandra's $c$-function first:
\proclaim{Definition 2.3} For $\l\in\fh^*$ and $k\in \K$ we write  
$$
\tilde c(\l,k)=\prod_{\a \in R_+}
{\G(\l(\a^\vee))
\over \G(\l(\a^\vee)+k_\a)}
\tag2.7
$$
and 
$$
c(\l,k)={\tilde c(\l,k)\over \tilde c(\rho(k),k)}
\tag2.8
$$
\endproclaim
\proclaim{Theorem 2.4} ([8]; the different sign in the argument of the
$c$-function is due to a change of sign in the definition of $\tilde c$.)
Assume that $k\in \K$ satisfies (with $\b\in R$ the highest short root) 
$$
\rho(k)(\b^\vee)+k_\b+1>0
\tag2.9
$$
(cf. (1.6)), that $\l$ is regular
and that $\l( \k^\vee) + 1
\not= 0\ (\forall \k \in Q \backslash \{ 0 \} )$. If $x \in \fa_-$ we define 
$$
F(\l,k;x) = \sum_{w \in W} c(- w \l, k) \Phi (w \l,k;x).
$$
This function extends to a holomorphic function in $(\l,k,z)$ in the 
domain where $\l\in\fh^*$, $k\in \K_c$ such that $\text{Re}(k)$ satisfies (2.9), 
and $z=x+iy$ with $x,\, y\in \fa$ such that $|\a(y)|<2\pi\ \forall\a\in R$. 
Moreover, $F(\l,k;z)$ is $W$-invariant as a function of both $\l$ and $z$, and 
satisfies $F(\l,k;0)=1$.
\endproclaim
\demo{Proof} Well known, cf. [8] or [17], Theorem 6.3. Notice that one avoids 
the zeros of the $\tilde c(\rho(k),k)$ because of (2.9).
\enddemo
Apart from the important asymptotic behaviour with respect to $z$ of 
$F(\l,k;z)$ given in 
Theorem 2.4 we want to have a good control over the growth behaviour 
in $\l\in\fh^*$ of $F(\l,k;x)$ when $x$ is confined to a compact subset of 
$\fa$. This is of crucial importance for the Paley-Wiener theorem. In 
[17] we proved a uniform estimate in both $x$ and in $\l$ but in this strong form 
the result will no longer be true if $k_\a<0$. Fortunately we can use the 
calculus of shift operators to prove a weaker statement that is just 
sufficient to prove the Paley-Wiener theorem.
\proclaim{Theorem 2.5} Assume that $k\in \K$ is such that $\tilde c(\rho(k),k)\not=0$.
Let $D\subset\fa$ be compact and let $p\in Sym(\fh)$. There exist constants  
$C\in \R_+$ and $M\in\bold N$ such that $\forall \l\in \fh^*$, $\forall x\in D$:
$$
|\partial(p)F(\l,k;x)|\leq C(1+|\l|)^Me^{\text{max}_w\text{Re}(w\l(x))}
\tag2.10
$$
\endproclaim
\demo{Proof}
By a slight variation of [17], Corollary 6.2, we see that this result is true 
when $k_\a\geq 0\ \forall\a\in R$. We now argue inductively, using the following 
formula from the calculus of hypergeometric shift operators (the formula 
is similar in nature to formula (5.1) of [17], but we omit details):
$$
b(k)F(\l,k;x)=\e_+\{\prod_{\a\in R_+}D_{\a^\vee}(k)(\Delta(x) F(\l,k+1;x))\}
\tag2.11
$$
where $\e_+$ denotes symmetrization over $W$ with respect to $x$, $k+1$ is the 
multiplicity function such that $(k+1)_\a=k_\a+1$, $\Delta$ denotes the 
Weyl denominator 
$$
\Delta(x)=
\prod_{\a\in R}(e^{\a(x)\over 2}-e^{-\a(x)\over 2})
$$
and $b(k)={\tilde c(\rho(k),k)\over\tilde c(\rho(k+1),k+1)}$, which is a polynomial 
in $k$, nonvanishing under the assumption that $\tilde c(\rho(k),k)\not=0$. 
We may assume that $D$ is $W$-invariant and convex. Now observe that 
$$
\eqalign{
D_\xi(k)f(x)=\partial(\xi)f(x)-&\rho(k)(\xi)f(x)-\cr
&\sum_{\a\in R_+}
k_\a{\a(x)\a(\xi)\over 1-e^{-\a(x)}}\int_{t=0}^1
\partial(\a^\vee)f(x-t\a(x)\a^\vee)dt\cr}
$$
From this expression we easily see that (2.10) holds for $F(\l,k,x)$, assuming that 
(2.10) holds for  $F(\l,k+1;x)$.
\enddemo
\beginsection{3. Wave packets, and the tempered spectrum}

{\sl In the rest of the paper we assume that the condition (1.6) holds, 
unless stated otherwise.}
Recall the definition of Paley-Wiener functions. 
\proclaim{Definition 3.1} Given $x \in \fa$, let $C_x$ denote the convex
hull of $Wx$. Define the support function $H_x$ on $\fa^*$ as
follows: $H_x(\l )=\text {sup}_{y\in C_x}\l (y)$. An entire function
$\phi$ on $\fh^*$ is said to have Paley-Wiener type $x$ if $\forall N \in
\bold Z_+ \exists C\in \bold R_+$ such that
$$
\vert \phi (\l )\vert \leq C(1+\vert \l \vert )^{-N}
\text{e}^{H_x(-\text {Re}\l )}
$$
The space of functions of Paley-Wiener type $x$ is denoted by $PW(x)$.
We also use the notation $PW=\cup_{x\in \fa}PW(x)$.
\endproclaim

We are going to study the following wave packet operator. 
If $\phi\in PW$ we take a point $p\in\fa^*_+$ such that 
$$
p(\a^\vee)+k_\a>0\ \forall\a\in R
\tag3.1
$$
and we define for $x\in\fa_-$:
$$
\Cal J(\phi)(x)=(2\pi i)^{-n}|W|^{-1}\int_{p+i\fa^*}\phi(\l)
\Phi(\l,k;x)
{d\l\over c(\l,k)}
\tag3.2
$$
We observe right away that this definition is independent of the choice of 
the point $p$ as long as $p$ satisfies the condition (3.1). First of all, 
the integrand is holomorphic in the domain ${\text{Re}}\l(\a^\vee)+k_\a>0
\ \forall\a\in R$. 
Because of Proposition 2.2 (2) we have a uniform bound on $\Phi(\l,k;x)$, so it
is allowed to apply Cauchy's theorem by the Paley-Wiener estimates on 
$\phi$ and the following simple fact for the reciprocal of 
the $c$-function: if $k_\a<0\ \forall\a\in R$ then $c(\l,k)^{-1}$ is 
bounded on regions of the form $\text{Re}\l(\a^\vee)+k_\a\geq\epsilon>0$.
 
Next, let us consider what happens when we move $p$ towards the origin 
and thus cross the poles of $c(\l,k)^{-1}$. This process was studied 
in a somewhat simpler context in [9]. Fortunately it is easy to reduce the 
situation to that of [9], so that we can use the results of that paper 
without difficulty. Let us first recall some important notions of [9], 
indispensable for the formulations of our results. 

For $L\subset \fa^*$ an affine subspace we put $R_L=\{\a\in R\mid
L(\a^\vee)=\text{constant}\}$. If $\fa_L^*=\text{span}(R_L)$ 
then it is clear that 
$R_L=R\cap
\fa^*_L$ is a parabolic root subsystem of $R$.
 An affine subspace $L\subset \fa^*$ is defined to be
residual (or more precisely $(\fa^*,R,k)$-residual) by induction on the 
codimension
of $L$. The space $\fa^*$ itself is by definition a residual subspace. 
The affine
subspace $L\subset \fa^*$ with positive codimension is called residual 
if there is
a residual subspace $M\subset \fa^*$ with $M\supset L$ and
$\text{dim}(M)=\text{dim}(L)+1$ such that 
$$
\#\{\a\in R_L\backslash R_M\mid L(\a^\vee)=k_\a\}\geq \#
\{\a\in R_L\backslash R_M\mid L(\a^\vee)=0\}+1
\tag3.3
$$
A residual point is also called a distinguished (or more precisely 
$(\fa^*,R,k)$
distinguished) point.
If $L\subset \fa^*$ is an affine subspace with
$\text{codim}(L)=
\text{rank}(R_L)$ then $L=c_L+V^L$ with $c_L$ the center of $L$
determined by $\{c_L\}=L\cap V_L$ and $V^L$ the orthogonal complement of $V_L$
in $V$. When $L$ is a residual subspace we define the 
tempered form $L^{temp}=c_L+iV^L\subset \fh^*$. There is a complete
classification of the  residual subspaces for each root system, cf. [9] 
section 4.

\proclaim{Lemma 3.2} For every residual subspace $L$, the center 
$c_L$ lies in the convex hull of the $W$-orbit of $\rho(k)$. In 
particular, $c_L$ satisfies $|c_L(\a^\vee)|<1\ \forall \a\in R$.  
\endproclaim
\demo{Proof} We may assume that $c_L$ is dominant. 
Clearly it is sufficient to prove the lemma for the 
situation $L=c_L$, i.e. the case when $L$ is a distinguished point.
Let $\b_1,\dots,\b_N$ be the positive roots such that 
$L(\b_i^\vee)=k_{\b_i}$. Using the classification of distinguished 
points of [9] we see that for every simple root $\a$ there exists 
an index $i$ such that $\b_i\geq\a$ {\sl and} $|\b_i|=|\a|$ (if $R$ is  
simply laced this is obvious; unfortunately I do not see a uniform 
proof in the other cases). Hence $L(\a^\vee)\leq -k_\a=-\rho(k)(\a^\vee)$ for 
every simple $\a\in R$, which is even a stronger statement than what 
we wanted to prove. The last statement of the lemma follows from 
Proposition 1.1
\enddemo  

The $c$-function has a nice decomposition in terms of the 
$c$-function of the Yang system which was introduced in 
[9]. Indeed, define 
$$
c_Y(\l,k)=\prod_{\a\in R_+}\frac{\a^\vee(\l)+k_\a}{\a^\vee(\l)}
\tag3.4
$$
and define $c^Y(\l,k)$ by the decomposition
$$
c(\l,k)=c_Y(\l,k)c^Y(\l,k).
\tag3.5
$$
We have 
\proclaim{Lemma 3.3} Let $U$ be a $W$-invariant open convex set 
that contains the convex hull of the $W$-orbit of $\rho(k)$ 
and which is itself contained in $\{\l\in\fa^*\mid\forall\a\in R:\ 
|\l(\a^\vee)|<1+k_\a\}$. 
Both $c^Y(\l,k)$ and $c^Y(\l,k)^{-1}$ are 
holomorphic in the tube $U+i\fa^*$, and $c^Y(\l,k)^{-1}$ is bounded
on this domain. 
Moreover, the product $c^Y(\l,k)c^Y(-\l,k)$ is $W$-invariant, and 
positive on every tempered subspace $L^{temp}$.
\endproclaim
\demo{Proof} We can choose $U$ as indicated because of Proposition 1.1, 
and because of Lemma 3.2 we know that $L^{temp}\subset U+i\fa^*$ 
for every  residual subspace $L$. The rest of the statements follow 
directly from the explicit formula 
$$
c^Y(\l,k)=\tilde c(\rho(k),k)^{-1}
\prod_{\a\in R_+}\frac{\G(\l(\a^\vee)+1)}{\G(\l(\a^\vee)+k_\a+1)}
\tag3.6
$$
The positivity is proved as in [9], Theorem 3.13. This is based 
on the fact that $-c_L\in W(R_L)c_L$, cf. [9], Theorem 3.10.
\enddemo

With these preliminaries in mind it is now clear that we can 
apply the residue calculus as formulated in [9], section 3, to 
the wave packet (3.2). Indeed, we can take $p$ in $U$ (see 
Lemma 3.3), still satisfying (3.1). 
By Lemma 3.2 and the fact that $U$ is convex and $W$-invariant
it is clear that the contour shifts needed in the calculation of 
residues all take place inside the tube $U+i\fa^*$.  
Therefore all the poles of 
the integrand one has to reckon with are those lying in 
the tube $U+i\fa^*$, and by Proposition 2.2 and Lemma 3.3 
these come only from $c_Y(\l,k)^{-1}$, the 
reciprocal of the $c$-function of the 
Yang system.  To describe the result, we introduce the following 
function supported on $L^{temp}$. If $\l\in L^{temp}$ write 
$$
f_L(\l,k)=\tilde c(\rho(k),k)^2\frac{\prod^\prime|\G(\l(\a^\vee)+k_\a)|}
{\prod^\prime|\G(\l(a^\vee))|}
\tag3.7
$$
where the $\prod^\prime$ denotes the product over all roots 
such that the argument of the corresponding gamma factor is not 
identically equal to $0$ on $L^{temp}$.    

With these definitions we have 
\proclaim{Theorem 3.4}Let $\omega_L$ denote the measure on $L^{temp}$ 
which is invariant under 
translations by elements of $iV^L$, normalized such that if $F$ is a 
fundamental domain of $V^L\cap 2\pi P$ 
($P$ is the weight lattice of $R$) then 
the volume of $c_L+iF$ is 1.
If $\phi\in PW$ and is $W$-invariant, and $x\in\fa_-$, then $\Cal J(\phi)(x)$
can also be written in the form 
$$
\Cal J(\phi)(x)=\int\phi(\l)F(\l,k;x)
d\nu(\l,k)
\tag3.8
$$
where $\nu=\sum_L\nu_L$ (the sum taken over all the residual 
subspaces), and $\nu_L$ 
is a measure supported on $L^{temp}$ of the form 
$$
\nu_L(\l,k)=\g_L(k)f_L(\l,k)\omega_L(\l)
\tag3.9
$$
Here $\g_L(k)$ is a nonnegative rational constant, depending 
on the configuration of hyperplanes of poles of $c_Y^{-1}$. It is 
$W$-invariant ($\g_L(k)=\g_{w(L)}(k)$) and invariant for scaling 
($\g_L(k)=\g_L(tk)$ if $t\in (0,1]$).
This constant is
hard to compute in general, but the following cases are easy: 
$\g_{\fa^*}=|W|^{-2}$, and if $c$ is a regular distinguished point 
defined by $\b_i^\vee(c)+k_{\b_i}=0,\,i=1,\dots,n,\,(\b_i\in R)$ then 
there are two possibilities: if $c$ is a strictly positive 
combination of the roots $\b_i$ then 
$\g_c=|W|^{-2}ind^{-1}$, where $ind$ is the index of the lattice 
spanned by the coroots $\b_i^\vee$ in $Q^\vee$. Otherwise $\g_c=0$. 
The function $f_L$ 
is smooth, nonnegative and bounded on $L^{temp}$. 
\endproclaim
\demo{Proof} Everything follows by applying the 
steps in the proof of [9], Theorem 3.18 to  
$$
\Cal J(\phi)(x)=(2\pi i)^{-n}|W|^{-1}\int_{p+i\fa^*}
\frac{\phi(\l)}{c^Y(\l,k)c^Y(-\l,k)}c^Y(-\l,k)\Phi(\l,k;x)
{d\l\over c_Y(\l,k)}
\tag3.10
$$
and the fact that (Theorem 2.4)
$$
F(\l,k;x)=\sum_{w\in W}c_Y(-w\l,k)c^Y(-w\l,k)\Phi(w\l,k;x)
\tag3.11
$$
We omit the details concerning the normalizations of measures. 
This is simply a matter of carefully comparing with the 
residue formulas in [9]. The statements about $f_L$ follow 
from Lemma 3.3 and formula (3.8), in combination with [9], formula
(3.2) (and the remarks following formula (3.2)).
\enddemo

Already at this point we can deduce very powerful conclusions 
about the leading terms of hypergeometric functions, by 
application of [9], Corollary 3.7 and Corollary 3.8 
to formula (3.11). We obtain:
\proclaim{Corollary 3.5} If $\l\in \text{supp}(\nu)$ then $F(\l,k)$ has 
a convergent asymptotic expansion on $\fa_-$ of the form 
$$
F(\l,k;x)=\sum_{\mu\in W\l}\sum_{\k\in Q_+}
p_{\mu,\k}(\l,k;x)e^{\{\rho(k)+\mu+\k\}(x)}
\tag3.12
$$
where the first sum is taken over $\mu\in W\l$ such that $\text{Re}(\mu)\geq 0$ 
(in the usual dominance ordering). 
Here the $p_{\mu,\k}$ are polynomials 
in $x$. If $\{\l\}=L$ is a distinguished point 
and $\nu_L(L)\not=0$ then this can be sharpened: the  
$\mu$ that are used in (3.12) have strictly positive real parts now. 
\endproclaim
\proclaim{Definition 3.6} The support of $\nu(\l,k)$ is called the 
tempered spectrum. 
If ${\l}=L$ is a distinguished point in the tempered spectrum then 
we call $\l$ a cuspidal point. 
\endproclaim

We can apply the techniques of Casselman and Mili\v ci\' c [3] in order 
to obtain estimates for $F(\l,k)$ on $\fa$ if $\l$ is in the tempered 
spectrum. 
\proclaim{Corollary 3.7} Let $\l$ be in the tempered spectrum, and choose  
$\tau$ in the closure of $\R_+ R_+$ such that $\tau\leq \text{Re}(\mu)$ for 
all $\mu$ that occur in (3.12). Then $\forall x\in\overline{\fa_-}$ we have 
$$
|F(\l,k;x)|\leq M(1+|x|)^me^{\{\tau +\rho(k)\}(x)}
\tag3.13
$$  
for suitable constants $M$ and $m$. In particular, if $\l$ is a cuspidal 
point then $F(\l,k)\in L_p(\fa,\delta(k,x)dx)$ for all $1\leq p<2+\epsilon$
for some $\epsilon>0$.
\endproclaim

\beginsection{4. The Paley-Wiener Theorems}

We are now in the position to prove the Paley-Wiener theorem. 
Given $f\in\cci(\fa)^W$ we define its Fourier transform $\F(f)$ by 
$$
\F(f)(\l)=\int_{\fa}f(x)F(-\l,k;x)\delta(k;x)dx
\tag4.1
$$
(Recall the normalization of the roots in formula (1.5).)
\proclaim{Theorem 4.1} If the support of $f\in\cci(\fa)^W$ is 
contained in $C_x$ (the convex hull of $Wx$) for some $x\in\fa$, 
then $\F(f)\in PW(x)$ (cf. Definition 3.1), and $\F(f)$ is 
$W$-invariant.
\endproclaim
\demo{Proof}  We use the fact that 
if $p$ is a real $W$-invariant polynomial on $\fa^*$ and 
$p^\vee$ denotes the polynomial $p^\vee(\l)=p(-\l)$ then 
$$
(D_pf,g)=(f,D_{p^\vee}g)
$$
where $(\cdot,\cdot)$ is the inner product (1.5), and $f$ and $g$ are 
in ${C}^\infty(\fa)^W$, at least one of them having compact 
support (cf. Lemma 7.8 of [17]). Hence, by Theorem 2.5, 
for all real $W$-invariant 
polynomials $p$ on $\fa^*$ there exist constants $C$ and $M$ such that  
$$
|p(\l)||\F(f)(\l)|\leq C(1+|\l|)^Me^{H_x(-\text{Re}\l)}
$$
with $M$ independent of $p$. This proves the result.
\enddemo
We now proceed with the converse statement for the 
wave packet operator $\Cal J$ discussed in the previous section. 
\proclaim{Theorem 4.2} If $\phi\in PW(x)$ for some $x\in \fa$, 
and $\phi$ is $W$-invariant, then $\Cal J(\phi)\in \cci(C_x)^W$. 
\endproclaim
\demo{Proof} First write $\Cal J(\phi)$ in the form (3.8). Using the 
boundedness of the functions $f_L$ on $L^{temp}$ and the estimates 
of Theorem 2.5 it is clear that $\Cal J(\phi)\in {C}^\infty(\fa)^W$. 
Next write,  
for $x\in\fa_-$, $\Cal J(\phi)(x)$ in the form (3.2). Now it is obvious that 
we may send $p$ to infinity in $\fa^*_+$ and apply the usual arguments 
of Helgason (cf. [11], Theorem 7.3) to see that the support of $\Cal J(\phi)$ 
is indeed contained in $C_x$. 
\enddemo

\beginsection{5. The inversion formula and the Plancherel formula}

We want to prove that $\Cal J$ is the inverse of $\F$. The proof  
is very similar to proof of this in the case when $k_\a>0$, as in 
[17], section 9. The heart of the matter is to use Peetre's 
characterization of differential operators [18].  The idea to use the 
result of Peetre goes back to Van den Ban and Schlichtkrull [1]. 
\proclaim{Lemma 5.1} The composition $\K=\Cal J\circ\F$ is 
formally symmetric on $\cci(\fa)^W$ with respect to $(\cdot,\cdot)$ 
(as defined 
in (1.5)). Moreover, $\K$ commutes with 
the operators $D_p$, $p\in Sym(\fh)^W$.  
\endproclaim
\demo{Proof} By the previous section $\K$ maps $\cci(\fa)^W$ onto itself. 
Using Theorem 3.4, Theorem 4.1 and Fubini's theorem we 
have for $f,g\in\cci(\fa)^W$:
$$
(\K f,g)=\int \F(f)(\l)\overline{\F(g)(\l)}d\nu(\l)
\tag5.1
$$
(we used that $-c_L\in W(R_L)c_L$, cf [9], Theorem 3.10.) 
and this is visibly a symmetric expression. 
The rest is trivial.
\enddemo
\proclaim{Proposition 5.2} The operator $\K$ is support preserving 
on $\cci(\fa)^W$. 
\endproclaim
\demo{Proof} By Lemma 5.1 and the previous section, $\K$ satisfies 
exactly the properties which are used in the argument of 
[17], Lemma 9.3.
\enddemo
\proclaim{Lemma 5.3} There exists a $W$-invariant polynomial $p$ on 
$\fa^*$ such that $\K=D_p$. 
\endproclaim
\demo{Proof} By Peetre's theorem [18] it is clear that $\K$ 
is a differential operator on $\cci(\fa^{\text{reg}})^W$, locally of 
finite order and with coefficients in $C^\infty(\fa^{\text{reg}})$. 
If one explicitly writes down the highest degree terms of 
the commutators $[D_{p_i},\K]$, ($\{p_i\}$ a set of generators of 
the ring of $W$-invariant polynomials) 
on a suitably small open set in $\fa^{\text{reg}}$ one sees that 
the highest degree term of $\K$ has to have constant coefficients. 
Because $\K$ extends to a map from $\cci(\fa)^W$ to itself one 
concludes that the highest order term of $\K$ is of the 
form $\partial(q)$ with $q$ a $W$-invariant polynomial. Now we can 
replace $\K$ by $\K-D_q$, which still commutes with all operators 
of the form $D_p$ but has lower order than $\K$. By induction 
on the order the result is proved.
\enddemo
\proclaim{Theorem 5.4} The maps $\F:\cci(\fa)^W\to PW^W$ and 
$\Cal J:PW^W\to \cci(\fa)^W$ are inverse to each other.
\endproclaim
\demo{Proof} To prove that $\K$ is the identity, we resort to the 
asymptotic formulas of Van den Ban and Schlichtkrull [1], 
Lemma 12.15 and Corollary 13.3 (as in [17], Lemma 9.8 and 
Lemma 9.9). The details are a little bit more complicated than 
in [17], 
because we have to reduce to the most continuous part of the 
spectrum first. Therefore, let $q$ be a nonzero $W$-invariant polynomial 
on $\fh^*$ that vanishes on the lower dimensional residual subspaces, 
and $p$ as in Lemma 5.3, then the argument of [17], Theorem 9.10 
applied to $\K\circ D_q$ shows that $pq=q$, hence $p=1$. The 
proof that $\K^\prime=\F\circ\Cal J$ is the identity is 
completely similar to [17], Lemma 9.11.
\enddemo
\proclaim{Theorem 5.5 (Plancherel Theorem)} The map $\F$ extends 
naturally to an isomorphism
$$
\F:L_2(\fa,\delta(k,x)dx)^W\to\left\{\bigoplus_{L\text{ tempered}}
L_2(L^{temp},\nu_L(k))\right\}^W.
\tag5.2
$$
\endproclaim
\demo{Proof} By (5.1) and Theorem 5.4 it is clear that $\F$ extends to 
an isometry defined on $L_2(\fa,\delta(k,x)dx)^W$. In the remaining part 
of the proof we will use $\F$ to denote this extended map.  
It remains to prove the surjectivity of $\F$. Note that, because $\F$ is  
an isometry defined on a Hilbert space, the image of $\F$ is closed. 
This being said, we embark on the proof of the surjectivity. First 
we claim that if $p$ is a nonvanishing polynomial on $i\fa^*$, then 
$p.PW$ resticted to $i\fa^*$ is dense in $L_2(i\fa^*,\nu_{i\fa^*})$. In 
order to see this, recall that $PW$ is dense in $\Cal S$, the ordinary 
Schwartz space on $i\fa^*$. If $\phi\in\cci(i\fa^*)$ and $f_n\to\phi$ in the 
Schwartz space topology, with $f_n\in PW$, then $p.f_n\to p.\phi$ in 
$L_2(i\fa^*,\nu_{i\fa^*})$ (because the function $c(\l,k)^{-1}$ is bounded 
on $i\fa^*$). On the other hand, by the nonvanishing of $p$, $p.\phi$ is 
an arbitrary element in $\cci(i\fa^*)$. Whence the density statement is 
proven, as claimed. 

There exists a $W$-invariant polynomial $p$ 
which is nonvanishing on $i\fa^*$ but which vanishes on all the residual subspaces 
other than $i\fa^*$ (see [9], proof of Lemma 3.1, uniqueness part). 
It is clear using Theorem 5.4 that 
the image of $\F$ contains $p.PW^W$ 
(where $p$ is the invariant polynomial just mentioned), 
and the closure of this in the right hand side of (5.2) equals $L_2(i\fa^*,\nu_{i\fa^*})^W$, 
by the above argument. Hence, the summand $L_2(i\fa^*,\nu_{i\fa^*})^W$ of the RHS 
of (5.2) is in the image of 
$\F$. Now we argue by induction on the length $|c_L|$ of a center $c_L$ of a  
residual subspace $L$. Assume that all the summands of the RHS of 
(5.2) that correspond to residual subspaces with center $c$ such that $|c|<|c_L|$ 
are in the image of $\F$. Let $p$ be a $W$-invariant polynomial which vanishes on 
all residual subspaces with center $c$ such that $|c|\geq |c_L|$ but $c\not\in Wc_L$, 
and nonvanishing on $W(c_L+i\fa^*)$ (again, such a polynomial exists by 
[9], proof of Lemma 3.1, uniqueness part). It follows 
that the restriction of $p.PW^W$ to 
the union of the tempered forms 
$L_1^{temp},...,L_l^{temp}$ of the residual subspaces contained in $W(c_L+i\fa^*)$, 
is in the image of $\F$. Call this space of functions $\Sigma_L$. 
Let $\phi$ be a $W$-invariant $\cci$ function 
on the $W$ orbit $W(c_L+i\fa^*)$, and $\phi^\prime$ its restriction to the 
union of $L_1^{temp},...,L_l^{temp}$. As before, $\phi$ can be approximated, with 
respect to the topology of the direct sum of the Schwartz spaces of functions on the 
$W$-translates of $c_L+i\fa^*$, by elements from $p.PW^W$ (restricted to 
$W(c_L+i\fa^*)$). This implies the approximation of $\phi^\prime$ by elements of 
$\Sigma_L$, because (3.7) and Theorem 3.4 show that the Plancherel measure 
$\nu_{L_i(k)}$ is given by a polynomially bounded measure which is smooth on 
$L_i^{temp}$. Finally, it is clear that the functions $\phi^\prime$ form a dense 
subspace of $\{\oplus L_2(L_i^{temp},\nu_{L_i(k)})\}^W$ because of the smoothness of 
$\nu_{L_i(k)}$ on $L_i^{temp}$. This completes the induction step. 
\enddemo
\proclaim{Remark 5.6} This argument also works in the context of the 
Lieb-McGuire system of particles studied in [9], cf. Remark 3.21 of that paper.
\endproclaim
\proclaim{Theorem 5.7} The $W$-invariant, square integrable (or 
cuspidal) eigenfunctions 
of $L(k)$ are precisely those of the form $F(\l,k)$ with 
$\l$ a cuspidal point. (The cuspidal points are classified in [9], 
section 4. There are a finite number of such points for each root system and 
each multiplicity function $k$ satisfying the conditions of Proposition 1.1.)
\endproclaim
\demo{Proof} First of all, by the classification of 
the dominant cuspidal points in [9], 
it is easy to check that these are ordered by the dominance ordering.  
In particular, the $L(k)$ eigenvalues are separated. We already know 
that the corresponding hypergeometric functions are square 
integrable by Corollary 3.7. It remains to show that 
these are the only $L_2$ eigenfunctions, but this is obvious from 
Theorem 5.5, since the measures $\nu_L$ are absolutely continuous with 
respect to the Lebesgue measure on $L^{temp}$. 
\enddemo
\proclaim{Corollary 5.8} Let $\l=L$ be a cuspidal point, for some 
multiplicity function $k_0\in \K$. There exists at least one linear 
parameter family of $(\fa^*,R,k)$ distinguished points 
$\l(k)$ such that $\l(k_0)=\l$ and such that $\l(k)$ is cuspidal 
in an open simplex $\Sigma$ containing $k_0$.   
In the simply laced case this simplex is always equal to  
$(-\frac{1}{d_n},0)$ (with $d_n$ the Coxeter number of $W$) and in the 
other cases $\Sigma$ is bounded by two lines through $0\in\K$ and 
the line $\rho(k)(\b^\vee)+k_\b+1=0$ (as in (1.6)). The function 
$F(\l(k),k;\cdot)$ is square integrable for all $k$ in $\Sigma+i\K$. 
Let $R_z=\{\a\in R\mid \l(k)(\a^\vee)=0\ \forall k\in \K$\} and 
let $R_p=\{\a\in R\mid \l(k)(\a^\vee)+k_\a=0\ \forall k\in \K$\}. Then 
$$
\int_\fa F(\l(k),k)^2\delta(k,x)dx=
c\frac{\prod_{\a\in R_+}\G(\rho(k)(\a^\vee)+k_\a)^2\prod
_{\a\in R\backslash R_z}\G(\l(k)(\a^\vee))}{\prod_{\a\in R_+}\G(\rho(k)(\a^\vee))^2
\prod_{\a\in R\backslash R_p}\G(\l(k)(\a^\vee)+k_\a)}
\tag5.3
$$
Here the constant $c$ is rational, and of the form 
$$
c=\pm\g_L^{-1}(k)|W\l(k)|^{-1}
$$
(cf. Theorem 3.4 for the definition of the 
rational number $\g_L(k)$) 
where $k\in\Sigma$ is a generic point. 
\endproclaim
\demo{Proof} The process of taking residues of  
$$
\int_{p+i\fa^*}\phi(\l)\frac{d\l}{c_Y(\l,k)c_Y(-\l,k)}
\tag5.4
$$
as described in [9] is continuous with respect to $k$ in the 
sense that the local contribution at $\l$ equals the sum of the 
limits of the local contributions at the cuspidal points $\l^\prime(k)$ which  
have the property that $\l^\prime(k)\to\l$ as $k\to k_0$ along a generic line. 
To see this, first note that (5.4) is continuous in $k$. 
Then one checks that the individual limits of the local contributions at 
the $\l^\prime(k)$ 
exist, by inequality (3.17) of [9]. Thus one may now use the uniqueness 
property of local contributions (cf. [9], Lemma (3.1)) to conclude the 
asserted continuity. Hence there exists at least one linear family 
$\l(k)$ of cuspidal points as stated in the theorem. The polynomials 
$p_{\mu,\k}(x;k)$ in the asymptotic  
expansion (3.12) of $F(\l(k),k,\cdot)$ are meromorphic in $k$, and their 
joint sets of singularities form an analytic subset 
of $\K_c$. Therefore, the set 
of leading terms of $F(\l(k),k,\cdot)$ 
at any specialization of $k\in \K_c$ is contained in the  
specialization of the generic set of exponents, and this generic set 
is of the form $\{w\l(k)\mid w\in D\}$ for some  
$D\subset W$. Moreover, on an arbitrary, sufficiently small 
compact torus $T$ in $\K_c$, with center $k_1$ say, which avoids 
the singularities of the $p_{\mu,\k}(x;k)$, the estimate 
(3.13) holds uniformly: $\forall k\in T\ \forall x\in
\overline{\fa_-}$: 
$$
|F(\l(k),k;x)|\leq M(1+|x|)^me^{\tau(x)}
$$
for some $\tau>\rho(k_1)$. Hence by Cauchy's theorem the same 
type of estimate holds for derivatives of $F(\l(k),k;x)$ with respect 
to $k$, and we conclude that the left hand side of (5.3) 
is a holomorphic function in in the tube 
$\{k\in \K_c\mid \forall w\in D: \text{Re}(w\l(k))>0\}\cap
\{k\in \K_c\mid \text{Re}(k)\text{ satisfies (1.6)}\}$. This set is of 
the form $\Sigma+i\K$ with $\Sigma$ a simplex as 
described in the theorem. 
We now compute the integral explicitly when $k\in\Sigma$ and generic,   
by applying $\F$ to $F(\l(k),k)$, 
and using Theorem 5.5. The rest of the statement 
follows by analytic continuation.
\enddemo

The following special case was previously obtained in [2], 
using different methods.
\proclaim{Example 5.9} Take $\l(k)=\rho(k)$ in the previous theorem. In this 
situation we know that  
$F(\rho(k),k;x)=1$, and thus 
$\Sigma=\{k\in \K\mid k\text{ satisfies (1.6)}\}$. The 
constant $c$ equals $|W|$, and if one takes $k_\a=k\ \forall \a\in R$ 
one obtains
$$
\int_{\fa}\delta(k,x)dx=
\prod_{i=1}^n\binom{d_ik}{k}\frac{\pi}{\sin{(-m_i\pi k)}}
$$
where the $d_i$ denote the primitive degrees of $W$, and the 
$m_i$ are the exponents of $W$. 
\endproclaim


\Refs

\ref
\no 1
\by Ban, van den, E.P., Schlichtkrull, H.
\paper The most continuous part of the Plancherel decomposition
		 for a reductive symmetric space
\jour Preprint series, University of Utrecht
\vol 921
\year 1995 
\endref
\ref 
\no 2
\by Brusse, R., Heckman, G.J., Opdam, E.M.
\paper Variation on a theme of Macdonald
\jour Math. Z.
\vol 208
\pages 1--10
\yr 1991
\endref
\ref
\no 3
\by Casselman, W., Mili\v ci\' c, D. 
\paper Asymptotic behaviour of matrix coefficients of admissible representations
\jour Duke Math. Jour. 
\vol 49
\pages 869--930
\yr 1982
\endref
\ref
\no 4
\by Cherednik, I.
\paper A unification of Knizhnik- Zamolodchikov equations 
    and Dunkl operators via affine Hecke algebras
\jour Inv. Math.
\vol 106
\pages 411--432
\yr 1991
\endref
\ref
\no 5
\by Dunkl, C.F.
\paper Differential-difference operators associated
		    to reflection groups
\jour Trans. AMS
\vol 311 No.1 
\pages
\yr 1989
\endref
\ref
\no 6
\by Harish-Chandra
\paper Spherical functions on a semisimple Lie group I
\jour Am. J. Math. 
\vol 80
\pages 241--310
\yr 1958
\endref
\ref 
\no 7
\by Heckman, G.J.
\paper An elementary approach to the hypergeometric
		   shift operators of Opdam
\jour Inv. Math.
\vol 103
\pages 341--350
\yr 1991
\endref
\ref 
\no 8
\by Heckman, G.J., Opdam, E.M.
\paper Root systems and hypergeometric functions I
\jour Comp. Math.
\vol 64
\pages 329--352
\yr 1987
\endref
\ref
\no 9
\by Heckman, G.J., Opdam, E.M.
\paper Yang's system of particles and Hecke algebras
\jour Annals of Mathematics
\vol 145
\pages 139--173
\yr 1997
\endref
\ref
\no 10
\by Heckman, G.J., Schlichtkrull, H.
\book Harmonic analysis and special functions on symmetric spaces, 
\bookinfo Perspectives in Mathematics 16
\publ Academic Press 
\yr 1995
\endref
\ref 
\no 11
\by Helgason, S.
\book Groups and geometric analysis
\publ Academic press 
\yr 1984
\endref
\ref
\no 12
\by Koornwinder, T.H.
\paper Orthogonal polynomials in two variables which are eigenfunctions of two 
algebraically independent partial differential operators I-IV
\jour Indag. Math. 
\vol 36
\pages 48--66 
\yr 1974
\moreref
\pages 358--381
\endref
\ref
\no 13
\by Macdonald, I.G.
\paper A formal identity for affine root systems
\jour preprint
\yr (1996)
\endref
\ref 
\no 14
\by Moser, J.
\paper Three integrable Hamiltonian systems connected with isospectral 
deformation
\jour Adv. Math.
\vol 16
\pages 197--220
\yr 1975
\endref
\ref
\no 15
\by Olshanetsky, M.A., Perelemov, A.M.
\paper Completely integrable Hamiltonian systems connected with 
semi\-simple Lie algebras 
\jour Inv. Math. 
\vol 37
\pages 93--108
\yr 1976
\endref
\ref
\no 16
\by Opdam, E.M.
\paper Some applications of hypergeometric shift
		   operators
\jour Inv. Math.
\vol 98
\pages 1--18
\yr 1989
\endref
\ref
\no 17
\by Opdam, E.M.
\paper Harmonic analysis for certain representations of 
                   graded Hecke algebras
\jour Acta Mathematica
\vol 175
\pages 75--121
\yr 1995
\endref
\ref
\no 18
\by Peetre, J.
\paper Rectification \`a l'article "Une caract\'erisation
		   abstraite des op\'era\-teurs diff\'erentiels"
\jour Math. Scand. 
\vol 8
\pages 116-120
\yr 1960
\endref
\endRefs
\enddocument